\documentclass[article]{amsart}
\usepackage{cases}
\usepackage{amsmath, amsfonts, pifont, amssymb}
\usepackage{verbatim}
\usepackage[bookmarks=true]{hyperref}
\usepackage{geometry}
\newtheorem{theorem}{Theorem}[section]
\newtheorem{lemma}[theorem]{Lemma}
\newtheorem{proposition}[theorem]{Proposition}

\newcommand{\beq}{\begin{equation}}
\newcommand{\eeq}{\end{equation}}
\newcommand{\beqq}{\begin{equation*}}
\newcommand{\eeqq}{\end{equation*}}

\theoremstyle{definition}
\newtheorem{definition}[theorem]{Definition}

\theoremstyle{remark}
\newtheorem{remark}[theorem]{Remark}

\numberwithin{equation}{section}

\begin{document}

\title[Dynamics for cubic NLS]{Long time dynamics for defocusing cubic NLS on three dimensional product space}
\author{Zehua Zhao and Jiqiang Zheng}
\maketitle

\begin{abstract}
In this article, we study long time dynamics for defocusing cubic NLS on three dimensional product space. First, we apply the decoupling method in Bourgain-Demeter \cite{BD} to establish a bilinear Strichartz estimate. Moreover, we prove global well-posedness for defocusing, cubic NLS on three dimensional product space with rough initial data ($H^s$, $s>\frac{5}{6}$) based on I-method and the bilinear estimate. At last, we discuss the growth of higher Sobolev norm problem which is tightly linked to the weak turbulence phenomenon.
\end{abstract}
\bigskip

\bigskip


\section{introduction}

We study the cubic, defocusing Schr\"odinger initial value problem on three dimensional product space as follows,
\begin{equation}\label{maineq}
\aligned
\begin{cases}
\left(i\partial_t+ \Delta_{\mathbb{R}^n \times \mathbb{T}^{3-n}}\right) u=  |u|^{2} u, \\
u(0,x,y) = u_{0}(x,y) \in H^{s}(\mathbb{R}^n \times \mathbb{T}^{3-n}),
\end{cases}
\endaligned
\end{equation}
where $n=0,1,2$, and $u(t,x,y):\mathbb{R}_t\times\mathbb{R}_x^n\times\mathbb{T}_y^{3-n}\to\mathbb{C}.$ Here the product space $\mathbb{R}^{m} \times \mathbb{T}^{n}$ is known as `semiperiodic space' as well as `waveguide manifold', where $\mathbb{T}^{n}$ is
a (rational or irrational) $n$-dimensional torus. This problem is known as low regularity problem when the Sobolev index satisfies $s<1$. Global well-posedness for \eqref{maineq} is well established according to classical contraction mapping method and conservation of energy if $s\geq 1$ (see \cite{Cbook,CW,Taobook} for examples), while delicate techniques are required for the low regularity case.\vspace{3mm}

According to the structure of semilinear NLS, there are three important conserved quantities of \eqref{maineq} as follows.
\begin{align*}
\text{Mass: }    &\quad
{M}(u(t))  = \int_{\mathbb{R}^n \times \mathbb{T}^{3-n}} |u(t,x,y)|^2\,\mathrm{d}x\mathrm{d}y,\\
\text{   Energy:  }     &  \quad
{E}(u(t))  = \int_{\mathbb{R}^n \times \mathbb{T}^{3-n}} \frac12 |\nabla u(t,x,y)|^2  + \frac{1}{4} |u(t,x,y)|^{4} \,\mathrm{d}x\mathrm{d}y,\\
\text{ Momentum: } &  \quad
{P}(u(t)) = \Im \int_{\mathbb{R}^n \times \mathbb{T}^{3-n}} \overline{u(t,x,y)} \nabla u(t,x,y)\,\mathrm{d}x\mathrm{d}y.
\end{align*}

In particular, we are interested in the long time dynamics of \eqref{maineq}. Generally, well-posedness theory and long time behavior of NLS is a hot topic in the area of dispersive evolution equations and has been studied widely in recent decades. Naturally, the Euclidean case is first treated and the theory at least in the defocusing setting has been well established. We refer to \cite{Iteam1,BD3,KM1} for some important  Euclidean results. Moreover, we refer to \cite{R2T,CGZ,CZZ,HP,HTT1,HTT2,IPT3,IPRT3,KV1,Haitian,Z1,Z2} with regard to the tori case and the waveguide case. One may roughly think that the  waveguide  case is  ‘between’  the Euclidean  case and the tori case in some sense since the waveguide is the product of the Euclidean space and the tori. Both of the techniques for the two cases are often combined and applied together to the waveguide case. \vspace{3mm}

One of the main theorems in this paper is as follows,
\begin{theorem}\label{main}
Initial value problem \eqref{maineq} is globally well-posed when $s>s_0=\frac{5}{6}$.
\end{theorem}

\begin{remark}
We believe this to be the first low regularity result for NLS on waveguide manifold. Additionally, we cover the 3d tori case as well (when $n=0$) in \eqref{maineq}. For 1d and 2d tori case, please see Silva-Pavlovic-Staffilani-Tzirakis \cite{Sgroup}.
\end{remark}

\begin{remark}
It's expected that Theorem \ref{main} holds for $s>\frac{1}{2}$ since \eqref{maineq} is ``$H^{\frac{1}{2}}$-critical''. Other methods or delicate techniques are required for one to obtain the sharp result.
\end{remark}

\begin{remark}
This result, (together with Theorem \ref{main2}), is indifferent to the rational/irrational choice for the tori direction since the tools and methods (such as Function spaces, Strichartz estimate and the bilinear estimate) we use work for both cases. For convenience of writing, in this paper, we discuss the normal case $\mathbb{T}^d=[0,1]^d$.
\end{remark}

\begin{remark}
For the purpose of unification, from now on, throughout this paper, we denote $\mathcal{M}$ to be the manifold $\mathbb{R}^n \times \mathbb{T}^{3-n}$ where $n=0,1,2$. When $n=3$, the NLS problem is on pure Euclidean space, which is quite different from other cases, so we dismiss this special case in our paper. Moreover, we denote $\mathcal{M}_{\lambda}$ to be the rescaled manifold $\mathbb{R}^n \times \mathbb{T}_{\lambda}^{3-n}$, where $\mathbb{T}_{\lambda}=[0,\lambda]$.
\end{remark}

The proof of Theorem \ref{main} follows from the classical I-method first established in Colliander-Keel-Staffilani-Takaoka-Tao \cite{Imethod} with a bilinear estimate which is based on decoupling method established in Bourgain-Demeter \cite{BD}. Please see \cite{Iteam3,BD1,Su,TR1,TR2} for some applications of I-method including some other more delicate techniques such as resonant decomposition, linear-nonlinear decomposition and potential bound control. We will study the bilinear estimate first (see Section 3) and then use it to prove decay of the modified energy (see Section 4). At last, we prove the low regularity result in Section 5.\vspace{3mm}

Moreover, we are interested in the growth of higher Sobolev norm problem which helps us understand the qualitative behavior of the solution better. This type of problems is hot in this area in recent decants and is related to the phenomenon of weak turbulence which is generally described as the solution transferring energy to higher and higher frequencies, causing the $H^s$ norm to grow while the $H^1$ norm remains bounded. Moreover, it is almost trivial to obtain an exponential upper bound for the $H^s$ norm, by iterating local in time theory. In Bourgain \cite{Bourgain2}, using his “high-low method”, Bourgain was first able to improve this to a bound that is polynomial in time, in the case of a cubic nonlinearity. In this paper, similar to the low regularity case, we apply an `upside-down' I-method instead in the setting of waveguide manifold. The following is our second main theorem.

\begin{theorem}\label{main2}
Suppose $u$ is a solution to \eqref{maineq} ($s=2$) with energy $E$ and $||u(0)||_{H^2}=A$. Then we have,
\begin{equation}
    ||u(t)||_{H^2}\lesssim A+(1+|t|)^{1+\delta},
\end{equation}
for any time $t$, and any $\delta>0$. Here all implicit constants will depend on $E$ and $\delta$, but not on $A$ or $t$.
\end{theorem}

\begin{remark}
The above result can be easily extended to the general $s>1$ case in the sense of replacing the exponent $1+\delta$ on the right hand side by $(s-1)+\delta$ without big changes. For the purpose of convenience, we prove the case when $s=2$. A three dimensional tori analogue is proved in Deng-Germain \cite{DG}.
\end{remark}

\begin{remark}
It is also interesting and meaningful to consider the lower bound problems for growth of higher Sobolev norms, which reveal the growth of some solutions by constructions. See Colliander-Keel-Staffilani-Takaoka-Tao \cite{Iteam2} for example.
\end{remark}

At last, we make a brief comment on this paper as well as our previous result Cheng-Zhao-Zheng \cite{CZZ}. Generally, we believe that for NLS problems, the waveguide case is `not worse' than the tori case in the following sense: the corresponding results in the waveguide setting would be the same as the tori case or even better. An evidence to view this is that the two basic estimates, i.e. the Strichartz estimate and the bilinear estimate in the waveguide setting is as same as in the tori setting. (See Barron \cite{Barron} and Section 3 of this paper respectively.) Thus, in our previous result \cite{CZZ}, we showed that we only need to care about the whole dimension of the waveguide, not the distribution of the Euclidean dimensions and the tori dimensions. Moreover, in this paper, we proved analogous results as the tori case for the three dimensional waveguide in two respects, i.e. low regularity result and growth of Sobolev norm result. One common point for the two problems is the I-method (a classical one and a `upside-down' one).  \vspace{3mm}

The next natural questions are `What is the difference between tori and waveguide (with the same whole dimension)?' and `How waveguide performs better than tori?'. There are some positive results to respond to these questions, which shows that the waveguide case is strictly `better' than the tori case for some situations. As an example, for the pure tori case, scattering behavior for NLS is not expected because of the lack of dispersion. However, it is possible in the waveguide setting for a proper choice of the geometry and the nonlinearity. See \cite{R2T,CGZ,HP,HTT1,HTT2,Z1,Z2} for examples.\vspace{3mm}

The organization of the rest of this paper is: in Section 2, we discuss the preliminaries including notations, function spaces and basic estimates. Moreover, we overview the setting of I-method; in Section 3, we establish the bilinear estimate in the setting of product space based on the decoupling method; in Section 4, we prove decay of the modified energy based on the bilinear estimate; in Section 5, we give the proof for Theorem \ref{main}; in Section 6, we discuss the growth of higher Sobolev norm problem and give proof for Theorem \ref{main2}.

\section{Preliminaries}
In this section, we will discuss notations, Strichartz estimate, Littlewood-Paley theory, function spaces and properties of I-operator and $\mathcal{D}$-operator.
\subsection{Notations and Definitions. }we write $A \lesssim B$ to say that there is a constant $C$ such that $A\leq CB$. We use $A \simeq B$ when $A \lesssim B \lesssim A $. Particularly, we write $A \lesssim_u B$ to express that $A\leq C(u)B$ for some constant $C(u)$ depending on $u$. In addition, $a\pm:=a\pm \epsilon$ with $0<\epsilon \ll 1$. We use $\chi$ to denote cutoff functions: compactly supported, and equal to one in a neighborhood of zero.\vspace{3mm}

Then we give some more preliminaries in the setting of waveguide manifold. The tori case can be defined similarly. In fact, it is included since it is a special case. Throughout this paper, we regularly refer to the spacetime norms
\begin{equation}
    ||u||_{L^p_tL^q_z(I_t \times \mathbb{R}^m\times \mathbb{T}^n)}=\left(\int_{I_t}\left(\int_{\mathbb{R}^m\times \mathbb{T}^n} |u(t,z)|^q dz \right)^{\frac{p}{q}} dt\right)^{\frac{1}{p}}.
\end{equation}
Moreover, we turn to the Fourier transformation and Littlewood-Paley theory. We define the Fourier transform on $\mathbb{R}^m \times \mathbb{T}^n$ as follows:
\begin{equation}
    (\mathcal{F} f)(\xi)= \int_{\mathbb{R}^m \times \mathbb{T}^n}f(z)e^{-iz\cdot \xi}dz,
\end{equation}
where $\xi=(\xi_1,\xi_2,...,\xi_{d})\in \mathbb{R}^m \times \mathbb{Z}^n$ and $d=m+n$. We also note the Fourier inversion formula
\begin{equation}
    f(z)=c \sum_{(\xi_{m+1},...,\xi_{d})\in \mathbb{Z}^n} \int_{(\xi_1,...,\xi_{m}) \in \mathbb{R}^m} (\mathcal{F} f)(\xi)e^{iz\cdot \xi}d\xi_1...d\xi_m.
\end{equation}
For convenience, we may consider the discrete sum to be the integral with discrete measure so we can combine the above integrals together and treat them to be one integral. Moreover, we define the Schr{\"o}dinger propagator $e^{it\Delta}$ by
\begin{equation}
    \left(\mathcal{F} e^{it\Delta}f\right)(\xi)=e^{-it|\xi|^2}(\mathcal{F} f)(\xi).
\end{equation}
We are now ready to define the Littlewood-Paley projections. First, we fix $\eta_1: \mathbb{R} \rightarrow [0,1]$, a smooth even function satisfying
\begin{equation}
    \eta_1(\xi) =
\begin{cases}
1, \ |\xi|\le 1,\\
0, \ |\xi|\ge 2,
\end{cases}
\end{equation}
and $N=2^j$ a dyadic integer. Let $\eta^d=\mathbb{R}^d\rightarrow [0,1]$, $\eta^d(\xi)=\eta_1(\xi_1)\eta_1(\xi_2)\eta_1(\xi_3)...\eta_1(\xi_d)$. We define the Littlewood-Paley projectors $P_{\leq N}$ and $P_{ N}$ by
\begin{equation}
    \mathcal{F} (P_{\leq N} f)(\xi):=\eta^d\left(\frac{\xi}{N}\right) \mathcal{F} (f)(\xi), \quad \xi \in \mathbb{R}^m \times \mathbb{Z}^n,
\end{equation}
and
\begin{equation}
P_Nf=P_{\leq N}f-P_{\leq \frac{N}{2}}f.
\end{equation}
For any $a\in (0,\infty)$, we define
\begin{equation}
    P_{\leq a}:=\sum_{N\leq a}P_N,\quad P_{> a}:=\sum_{N>a}P_N.
\end{equation}

Scaling is an important symmetry for NLS on Euclidean space. However, in the tori or waveguide setting, when doing a rescaling, one needs to be careful about the spatial domain which will be changed (extended or contracted). If $u$ solves \eqref{maineq} on $\mathbb{R}^n \times \mathbb{T}^{3-n}$, then $u^{\lambda}(t,x,y)=\frac{1}{\lambda}u(\frac{t}{\lambda^2},\frac{x}{\lambda},\frac{y}{\lambda})$ solves the equation on $\mathbb{R}^n \times \mathbb{T}_{\lambda}^{3-n}$ with $\mathbb{T}_{\lambda}=[0,\lambda]$.\vspace{3mm}

We write $U_\lambda(t)$ for the solution operator to the rescaled linear Schr\"odinger equation,
\begin{equation}
  \left(i\partial_t+ \Delta_{\mathbb{R}^n \times \mathbb{T}_{\lambda}^{3-n}}\right) u(t,z)=0, \textmd{  where,  }z\in \mathbb{R}^n \times \mathbb{T}_{\lambda}^{3-n}.
\end{equation}
Equivalently,
\begin{equation}
   U_\lambda(t)u_0(z)=\int e^{2\pi ik\cdot z-(2\pi k)^2it}\mathcal{F}(u_0)(k)(dk)_{\lambda},
\end{equation}
where $(dk)_{\lambda}$ is the measure corresponding to the rescaled Fourier space.
\subsection{Function spaces and basic estimates}

We recall the Fourier restriction space (also known as `Bourgain space') as follows,
\begin{equation}
    ||u||_{X^{s,b}}=||U_{\lambda}(-t)u||_{H^b_tH^s_x}=||\langle k \rangle^s \langle \tau-4\pi^2 k^2 \rangle^{b}\mathcal{F}(u)(k,\tau) ||_{L^2_{\tau}L^2_{(dk)_{\lambda}}}.
\end{equation}
For the sake of $L^{\infty}$-embedding and the transfer principle, $b$ is chosen to be $\frac{1}{2}+$ throughout this paper. See Tao \cite{Taobook} for more information and properties regarding Bourgain space.
\begin{lemma}[Strichartz estimate]\label{Stricharz}
We recall the Strichartz estimate for $\mathcal{M}_{\lambda}$ in the following form (see Barron \cite{Barron} for the waveguide case and see for Killip-Visan \cite{KV1} the tori case),
\begin{equation}\label{Stri}
    ||u||_{L^{\frac{10}{3}}_{t,x}}\lesssim \lambda^{0+}||u||_{X^{0+,\frac{1}{2}+}}.
\end{equation}
\end{lemma}
Similar to the 2d case (see \cite{Sgroup}), we have
\begin{equation}\label{infty}
    ||u||_{L^{\infty}_{t,x}}\lesssim ||u||_{X^{\frac{3}{2}+,\frac{1}{2}+}}.
\end{equation}
By interpolation, we obtain
\begin{equation}\label{interpolate}
    ||u||_{L^{p}_{t,x}}\lesssim \lambda^{0+}||u||_{X^{\alpha(p),\frac{1}{2}+}},
\end{equation}
where $\alpha(p)=(1-\frac{10}{3p})\frac{3}{2}+$. In particular,
\begin{equation}\label{interpolate1}
    ||u||_{L^{\frac{30}{7}}_{t,x}}\lesssim \lambda^{0+}||u||_{X^{\frac{1}{3}+,\frac{1}{2}+}},
\end{equation}
and
\begin{equation}\label{interpolate2}
    ||u||_{L^{\frac{15}{2}}_{t,x}}\lesssim \lambda^{0+}||u||_{X^{\frac{5}{6}+,\frac{1}{2}+}}\end{equation}
holds.
\subsection{Setting of `usual' I-method}
We define the I-operator as in previous low regularity results (see Colliander-Keel-Staffilani-Takaoka-Tao \cite{Imethod} as an example). Precisely, given a frequency parameter $N \gg 1$ to be chosen later, we define $m(k)$
to be a smooth and decreasing multiplier satisfying,
\begin{equation}
    m(k) =
\begin{cases}
1, k < N,\\
(\frac{N}{k})^{1-s}, k > 2N.
\end{cases}
\end{equation}
Then $I:H^s \rightarrow H^1$ will be defined as the following multiplier operator
\begin{equation}
    \widehat{Iu}(k)=m(k)\hat{u}(k).
\end{equation}
This operator is smoothing of order $1-s$, that is:
\begin{equation}\label{Ioper}
  ||u||_{X^{s_0,b_0}}  \lesssim ||Iu||_{X^{s_0+1-s,b_0}} \lesssim N^{1-s}||u||_{X^{s_0,b_0}},
\end{equation}
for any $s_0,b_0\in \mathbb{R}$.
\begin{definition}[Modified energy]
    We define the modified energy as follows,
    \begin{equation}
        E(Iu)=\frac{1}{2}\int |\nabla I u|^2+\frac{1}{4}\int |Iu|^4dx.
    \end{equation}
\end{definition}
There is a basic estimate about the symbol $m(\xi)$. For $k \gtrsim N$ and $\alpha \geq (1-s)-$,
\begin{equation}\label{Ibasic}
    \frac{1}{m(k)|k|^{\alpha}}\lesssim N^{-\alpha}.
\end{equation}
A key step is to investigate the `almost' conserved properties of the modified energy which will be explicitly discussed in Section $4$.

\subsection{Setting of `upside down' I-method}
Similar to the idea of I-operator, in order to deal with higher Sobolev norm, we define the $\mathcal{D}$-operator as follows. For any fixed scale $N$, let the multiplier $\mathcal{D}$ to be
\begin{equation}
    \widehat{\mathcal{D}u}(k)=m(k)\hat{u}(k),
\end{equation}
where $m(k)$ is a smooth, even and decreasing function satisfying
\begin{equation}
    m(k) =
\begin{cases}
1, k < N,\\
\frac{k}{N}, k > 2N.
\end{cases}
\end{equation}

\section{Bilinear Strichartz estimate}
In this section, we establish a bilinear Strichartz estimate in the waveguide setting based on decoupling method in Bourgain-Demeter \cite{BD}, which is a crucial step to obtain the low regularity result. First, we give a brief overview of existing related results. For convenience, we recall them in the setting of three dimensional case and we refer to corresponding original manuscripts for general cases.\vspace{3mm}

First, for Euclidean case, one has,
\begin{proposition}[Bilinear estimate in the Euclidean setting]\label{BiEuc}
Fix $0<T<1$, then for $1\leq N_2 \leq N_1$, we have
\begin{equation}\label{bilinear1}
    ||e^{it\Delta}u_{N_1} e^{it\Delta}v_{N_2}||_{L^2_{t,x}([0,T)\times \mathbb{R}^{3} )}\lesssim N_2N_1^{-\frac{1}{2}}||u_{N_1}||_{L^2} ||v_{N_2}||_{L^2}.
\end{equation}
\end{proposition}
The above estimate is classical and has widely been applied. See Bourgain \cite{Bourgain1998} for instance. For tori case, there are two results as follows.
\begin{proposition}[Bilinear estimate in Killip-Visan \cite{KV1}]\label{Bi1}
Fix $0<T<1$, then for $1\leq N_2 \leq N_1$, we have
\begin{equation}\label{bilinear1}
    ||e^{it\Delta}u_{N_1} e^{it\Delta}v_{N_2}||_{L^2_{t,x}([0,T)\times \mathbb{T}^{3} )}\lesssim N_2^{\frac{1}{2}}||u_{N_1}||_{L^2} ||v_{N_2}||_{L^2}.
\end{equation}
This implies,
\begin{equation}
    ||u_{N_1} v_{N_2}||_{L^2_{t,x}([0,T)\times  \mathbb{T}^{3} )}\lesssim N_2^{\frac{1}{2}}||u_{N_1}||_{X^{0,1/2+}} ||v_{N_2}||_{X^{0,1/2+}}.
\end{equation}
\end{proposition}

\begin{proposition}[Bilinear estimate in Fan-Staffilani-Wang-Wilson \cite{Fan}]\label{Bi2}
Fix $0<T<1$, then for $1\leq N_2 \leq N_1$, we have
\begin{equation}
    ||U_\lambda(t)u_{N_1} U_\lambda(t)v_{N_2}||_{L^2_{t,x}([0,T)\times \mathbb{T}_{\lambda}^{3} )}\lesssim N_2^{\epsilon}\Big(\frac{1}{\lambda^{\frac{1}{2}}}+\frac{N_2}{N_1^{\frac{1}{2}}}\Big)||u_{N_1}||_{L^2} ||v_{N_2}||_{L^2}.
\end{equation}
This implies,
\begin{equation}\label{bilinear4}
    ||u_{N_1} v_{N_2}||_{L^2_{t,x}([0,T)\times  \mathbb{T}_{\lambda}^{3} )}\lesssim N_2^{\epsilon}\Big(\frac{1}{\lambda^{\frac{1}{2}}}+\frac{N_2}{N_1^{\frac{1}{2}}}\Big)||u_{N_1}||_{X^{0,1/2+}} ||v_{N_2}||_{X^{0,1/2+}}.
\end{equation}
\end{proposition}

\begin{remark}
Comparing Proposition \ref{BiEuc} with Proposition \ref{Bi1}, it is clear that in the Euclidean setting, one has stronger estimate. Moreover, comparing Proposition \ref{Bi2} with Proposition \ref{BiEuc}, letting $\lambda$ goes to infinity and dismissing the loss of $\epsilon$, the estimates are consistent with the bilinear Strichartz inequality in the Euclidean setting.
\end{remark}

\begin{remark}
In this squeal for $\mathbb{T}^3$ case, it is better to apply Proposition \ref{Bi2} instead of Proposition \ref{Bi1}, which will give us better result, more precisely, better decay of modified energy. Moreover, when we apply \eqref{bilinear4} to estimate terms in Proposition \ref{almost}, $\frac{1}{\lambda^{\frac{1}{2}}}$-term is comparably much easier to be handled. So it suffices to focus on the second term.
\end{remark}

\begin{remark}
In Silva-Pavlovic-Staffilani-Tzirakis \cite{Sgroup}, there are bilinear estimates as well and they are proved by using number theory techniques. Recent decoupling method established in Bourgain-Demeter \cite{BD} allows people to derive bilinear estimate as in Fan-Staffilani-Wang-Wilson \cite{Fan}.
\end{remark}

Thus now the next task is to generalize Proposition \ref{Bi2} for product space case which would be enough for manifold $\mathcal{M}$ concerning this paper since the 3d tori case has already been covered by Proposition \ref{Bi2}. For product space, one has
\begin{proposition}[Bilinear estimate in Cheng-Zhao-Zheng \cite{CZZ}]\label{Bi3}
Fix $0<T<1$, then for $1\leq N_2 \leq N_1$, we have
\begin{equation}\label{bilinear3}
    ||e^{it\Delta}u_{N_1} e^{it\Delta}v_{N_2}||_{L^2_{t,x}([0,T)\times \mathcal{M} )}\lesssim N_2^{\frac{1}{2}}||u_{N_1}||_{L^2} ||v_{N_2}||_{L^2}.
\end{equation}
This implies,
\begin{equation}
    ||u_{N_1} v_{N_2}||_{L^2_{t,x}([0,T)\times  \mathcal{M} )}\lesssim N_2^{\frac{1}{2}}||u_{N_1}||_{X^{0,1/2+}} ||v_{N_2}||_{X^{0,1/2+}}.
\end{equation}
\end{proposition}
The above estimate is the waveguide analogue of Proposition \ref{Bi1}. However, it is not as strong as the analogue of Proposition \ref{Bi2} in some sense. We expect to have the following estimate,
\begin{proposition}[Bilinear estimate]\label{Bi4}
Fix $0<T<1$, then for $1\leq N_2 \leq N_1$, we have
\begin{equation}
    ||U_\lambda(t)u_{N_1} U_\lambda(t)v_{N_2}||_{L^2_{t,x}([0,T)\times \mathcal{M}_{\lambda} )}\lesssim N_2^{\epsilon}\Big(\frac{1}{\lambda^{\frac{1}{2}}}+\frac{N_2}{N_1^{\frac{1}{2}}}\Big)||u_{N_1}||_{L^2} ||v_{N_2}||_{L^2}.
\end{equation}
This implies,
\begin{equation}\label{bilinear2}
    ||u_{N_1} v_{N_2}||_{L^2_{t,x}([0,T)\times  \mathcal{M}_{\lambda} )}\lesssim N_2^{\epsilon}\Big(\frac{1}{\lambda^{\frac{1}{2}}}+\frac{N_2}{N_1^{\frac{1}{2}}}\Big)||u_{N_1}||_{X^{0,1/2+}} ||v_{N_2}||_{X^{0,1/2+}}.
\end{equation}
\end{proposition}
The proof of the above estimate is based on decoupling method and has similar spirit of Fan-Staffilani-Wang-Wilson \cite{Fan} with some minor modifications. In fact, the $\mathbb{T}^3$ case is already included and we only need to take care of the waveguide case. We aim to prove the following more general result.
\begin{proposition}[Bilinear estimate in the setting of waveguide manifold]\label{Bi5}
Fix $0<T<1$, then for $1\leq N_2 \leq N_1$, we have
\begin{equation}\label{Bimain}
    ||U_\lambda(t)u_{N_1} U_\lambda(t)v_{N_2}||_{L^2_{t,x}([0,T)\times \mathbb{R}^n \times \mathbb{T}_{\lambda}^{d-n} )}\lesssim N_2^{\epsilon}D_{\lambda,N_1,N_2}||u_{N_1}||_{L^2} ||v_{N_2}||_{L^2},
\end{equation}
where (for convenience of notation)
\begin{equation}
D_{\lambda,N_1,N_2}:=
    \begin{cases}
(1/\lambda+N_2/N_1)^{\frac{1}{2}}, d=2,\\
(N_2^{d-3}/\lambda+N_2^{d-1}/N_1)^{\frac{1}{2}}, d > 3.
\end{cases}
\end{equation}
\end{proposition}
We will give an overview of the proof for Proposition \ref{Bi5} as follows. The main idea is same as the pure tori case shown in Fan-Staffilani-Wang-Wilson \cite{Fan}. We start by a lemma which has the same spirit as the standard parallel decoupling in the bilinear setting.
\begin{lemma}\label{paral}
We divide the domain $Q:=[0,T)\times \mathbb{R}^n \times \mathbb{T}_{\lambda}^{d-n} $ into the disjoint union of $Q_l:=[0,T)\times B_l \times \mathbb{T}_{\lambda}^{d-n}$ where
$B_l$'s are cubes with length $\lambda $ in $\mathbb{R}^n$. In order to prove Proposition \ref{Bi5} in the setting of $Q$, it suffices to prove it for each $Q_l$ in the following sense,
\begin{equation}\label{Bimain2}
    ||U_\lambda(t)u_{N_1} U_\lambda(t)v_{N_2}||_{L^2_{t,x}([0,T)\times B_l \times \mathbb{T}_{\lambda}^{d-n} )}\lesssim N_2^{\epsilon}D_{\lambda,N_1,N_2}||u_{N_1}||_{L^2(B_l \times \mathbb{T}_{\lambda}^{d-n})} ||v_{N_2}||_{L^2(B_l \times \mathbb{T}_{\lambda}^{d-n})}.
\end{equation}
\end{lemma}
\emph{Proof of Lemma \ref{paral}: }The proof is straightforward which is based on the decomposition and Minkowski inequality. Now we assume that \eqref{Bimain2} holds for all $B_l$. Then, by decomposition and \eqref{Bimain2},
\begin{equation}
    \aligned
        ||U_\lambda(t)u_{N_1} U_\lambda(t)v_{N_2}||^{2}_{L^2_{t,x}([0,T)\times \mathbb{R}^n \times \mathbb{T}_{\lambda}^{d-n} )}& \lesssim \sum_{l}         ||U_\lambda(t)u_{N_1} U_\lambda(t)v_{N_2}||^{2}_{L^2_{t,x}([0,T)\times B_l \times \mathbb{T}_{\lambda}^{d-n} )} \\
        &\lesssim \sum_{l} N_2^{2\epsilon}D^2_{\lambda,N_1,N_2}||u_{N_1}||^2_{L^2(B_l \times \mathbb{T}_{\lambda}^{d-n})} ||v_{N_2}||^2_{L^2(B_l \times \mathbb{T}_{\lambda}^{d-n})}           \\
        &\lesssim N_2^{2\epsilon}D^2_{\lambda,N_1,N_2} ||u_{N_1}||^2_{L^2} ||v_{N_2}||^2_{L^2}.
    \endaligned
\end{equation}
This implies the estimate \eqref{Bimain} in Proposition \ref{Bi5}.\vspace{3mm}

Now we recall Theorem 1.3 in Fan-Staffilani-Wang-Wilson \cite{Fan} as follows.
\begin{theorem}[Bilinear Decoupling Argument]\label{BilDec}
Let $f_1$ be supported on $P$ where $|\xi| \sim 1$, and let $f_2$ be
supported on $P$ where $|\xi| \sim N_2/N_1$. Let domain $\Omega=\{(t,x)\in [0,N_1^2]\times [0,(\lambda N_1)^2]^d\}$. For a finitely overlapping covering of the ball $B=\{|\xi| \leq 1\}$ of caps $\{ \theta \}$, we have the following estimate.
\begin{equation}\label{Bimain3}
    ||Ef_1 Ef_2||_{L_{avg}^2(\omega_{\Omega})}\lesssim N_2^{\epsilon}\lambda^{\frac{d}{2}} D_{\lambda,N_1,N_2}\prod_{j=1}^{2}(\sum_{|\theta|=1/(\lambda N_1)}||Ef_{j,\theta}||^2_{L^4_{avg}(\omega_{\Omega})})^{\frac{1}{2}}.
\end{equation}
\end{theorem}
As for the pure tori case, the above theorem is essential to the bilinear estimate and most of the paper \cite{Fan} details the proof of Theorem \ref{BilDec}. More precisely, Proposition \ref{Bi4} is implied by Theorem \ref{BilDec} by a standard reduction process as in Bourgain-Demeter \cite{BD}. We refer to \cite{Fan} for explanation of the notations, background, more details and motivations about this method.\vspace{3mm}

As for the waveguide case, our bilinear result Proposition \ref{Bi5} will be obtained by combining Theorem \ref{BilDec} and Lemma \ref{paral} using a similar reduction process. The strategy is to reduce the waveguide case to the tori case and then we can follow the standard process. We refer to \cite{Fan} for more details. The idea of the proof are stated as follows.\vspace{3mm}

\emph{Proof of Proposition \ref{Bi5}: }It suffices to prove \eqref{Bimain2} according to Lemma \ref{paral}. Without loss of generality, we set the cube $B_l:=[0,\lambda ]^n$. Then, similar to the tori case, to make the settings in Proposition \ref{Bi5} and Theorem \ref{BilDec} compatible, we rescale $\phi_1$ to be supported in the unit ball and rescale $\phi_2$ to be
supported in a ball of radius $\sim \frac{N2}{N1}$.\vspace{3mm}

Now we consider the set $Q_0=[0,N_1^2]\times [0,\lambda N_1]^n \times \mathbb{T}^{d-n}_{\lambda N_1}$. Here, we view $\mathbb{T}^{d-n}_{\lambda N_1}$ as a compact set in $\mathbb{R}^{d-n}$. So we may treat $[0,\lambda N_1]^n$ and $\mathbb{T}^{d-n}_{\lambda N_1}$ together without any differences, which implies the reduction to the tori case. Moreover, it is still noted that $Q_0$ is obviously smaller than $\Omega$ in Theorem \ref{BilDec} and $\Omega$ can be covered by $Q$ such that $\{ Q \}$ are finitely overlapping and each $Q$ is a translation of $Q_0$. \vspace{3mm}

Then it is clear that the rest of the proof follows exactly as the tori case so we omit it. See Fan-Staffilani-Wang-Wilson \cite{Fan} for more details.

\section{Decay of the Modified Energy}
Now we are ready to discuss the decay of the modified energy. H\"older's inequality, Strichartz estimate, estimate \eqref{Ibasic} and Bilinear estimate \eqref{bilinear2} are frequently used.
\begin{proposition}[Almost conservation law]\label{almost}
Let $u$ be a solution to \eqref{maineq}, then
\begin{equation}
    |E(Iu)(t)-E(Iu)(0)| \lesssim \frac{1}{N^{1-}}||Iu||_{X^{0,\frac{1}{2}+}}.
\end{equation}
\end{proposition}
\begin{remark}
See Proposition 4.7 of Silva-Pavlovic-Staffilani-Tzirakis \cite{Sgroup} for the 2d tori analogue of this proposition.
\end{remark}
\emph{Proof of Proposition \ref{almost}: }Direct calculations imply,

\begin{equation}
    \partial_t E(Iu)(t)=\textmd{Re}\int_{\mathcal{M}_{\lambda}}\overline{Iu_t}(|Iu|^2Iu-I(|u|^2u)).
\end{equation}
Then by Parseval’s formula, integrate in time, we have
\begin{equation}
    E(Iu)(t)-E(Iu)(0)=\int_0^t \int_{k_1+k_2+k_3+k_4=0} \left( 1-\frac{m(k_2+k_3+k_4)}{m(k_2)m(k_3)m(k_4)} \right) \widehat{\overline{I\partial_t u}}(k_1) \widehat{{I u}}(k_2) \widehat{\overline{I u}}(k_3) \widehat{\overline{I u}}(k_4).
\end{equation}
Using equation \eqref{maineq}, we split the above into two terms and estimate them separately. We let

\begin{equation}\label{term1}
    I:=\int_0^t \int_{k_1+k_2+k_3+k_4=0} \left( 1-\frac{m(k_2+k_3+k_4)}{m(k_2)m(k_3)m(k_4)} \right) \widehat{\Delta \overline{ I u}}(k_1) \widehat{{I u}}(k_2) \widehat{\overline{I u}}(k_3) \widehat{\overline{I u}}(k_4),
\end{equation}
and
\begin{equation}\label{term2}
    II:=\int_0^t \int_{k_1+k_2+k_3+k_4=0} \left( 1-\frac{m(k_2+k_3+k_4)}{m(k_2)m(k_3)m(k_4)} \right) \widehat{\overline{I (|u|^2 u)}}(k_1) \widehat{{I u}}(k_2) \widehat{\overline{I u}}(k_3) \widehat{\overline{I u}}(k_4).
\end{equation}
We start by estimating the first term. First notice that,
\begin{equation}
    ||\Delta(Iu) ||_{X^{-1,\frac{1}{2}+}} \leq ||Iu||_{X^{1,\frac{1}{2}+}},
\end{equation}
thus, after breaking the functions into by Littlewood Pelay theory, it suffices to show,
\begin{equation}\label{bound1}
    \aligned
    & \int_0^t \int_{\Gamma_4} \left( 1-\frac{m(N_2+N_3+N_4)}{m(N_2)m(N_3)m(N_4)} \right) \widehat{ \overline{ \phi_1}}(k_1) \widehat{{\phi_2}}(k_2) \widehat{\overline{\phi_3}}(k_3) \widehat{\overline{\phi_4}}(k_4)                             \\
    &\lesssim \frac{1}{N^{1-}}(N_1N_2N_3N_4)^{0-} ||\phi_1||_{X^{-1,\frac{1}{2}+}}\prod_{i=2}^{4}||\phi_i||_{X^{1,\frac{1}{2}+}},
    \endaligned
\end{equation}
where $\phi_i$ ($i=1,2,3,4$) are $\lambda$-periodic function with positive spatial Fourier transforms supported on annulus with radius $N_i$.\vspace{3mm}

Without loss of generality, we can assume that $N_2 \leq N_3 \leq N_4$. This implies $N_1 \lesssim N_2$ noticing the hyperbolic plane relation. For the sake of analyzing the symbol, we will consider several scenarios by discussing the size of $N$.\vspace{3mm}

\emph{Case I: $N \gg N_2$. } This case is trivial since the symbol is identically zero which implies the bound.\vspace{3mm}

\emph{Case II: $N_2 \gtrsim N \gg N_3 \geq N_4 $. }
Noticing the hyperbolic plane relation, we have $N_1\sim N_2$. Also, the symbol in \eqref{bound1} has bound $\frac{N_3}{N_2}$. Using the pointwise bound for the symbol, the Cauchy-Schwartz inequality, and Plancharel’s theorem we obtain,
\begin{equation}
    \textmd{LHS of \eqref{bound1}}\lesssim \frac{N_3}{N_2} ||\phi_1 \phi_3||_{L^2_{t,x}} ||\phi_2 \phi_4||_{L^2_{t,x}}.
\end{equation}
By Bilinear estimate \eqref{bilinear2},
\begin{equation}
\aligned
    \textmd{LHS of \eqref{bound1}}&\lesssim \frac{N_3}{N_2} \frac{N_3N_4}{N_1^{\frac{1}{2}}N_2^{\frac{1}{2}}}N_3^{+}N_4^{+}\frac{N_1}{N_2 N_3 N_4}||\phi_1||_{X^{-1,\frac{1}{2}+}}\prod_{i=2}^{4}||\phi_i||_{X^{1,\frac{1}{2}+}} \\
    &\lesssim    \frac{1}{N^{1-}} N_2^{0-} ||\phi_1||_{X^{-1,\frac{1}{2}+}}\prod_{i=2}^{4}||\phi_i||_{X^{1,\frac{1}{2}+}}     .
    \endaligned
\end{equation}
So this case can be handled.\vspace{3mm}

\emph{Case III: $N_2 \geq N_3 \gtrsim N$. }We will use the following bound on the multiplier,
\begin{equation}\label{pointwise3}
   \left( 1-\frac{m(N_2+N_3+N_4)}{m(N_2)m(N_3)m(N_4)} \right)  \lesssim \frac{m(N_1)}{m(N_2)m(N_3)m(N_4)}.
\end{equation}
Using the pointwise bound \eqref{pointwise3}, the Cauchy-Schwartz inequality, and Plancharel’s theorem we obtain,
\begin{equation}
    \textmd{LHS of \eqref{bound1}}\lesssim \frac{m(N_1)}{m(N_2)m(N_3)m(N_4)} ||\phi_1 \phi_3||_{L^2_{t,x}} ||\phi_2 \phi_4||_{L^2_{t,x}}.
\end{equation}
Now we distinguish two subcases.\vspace{3mm}

\emph{Case III (a): $N_2 \sim N_3 \gtrsim N$. }
Moreover, by Bilinear estimate \eqref{bilinear2},
\begin{equation}
\aligned
    \textmd{LHS of \eqref{bound1}}&\lesssim \frac{m(N_1)}{m(N_2)m(N_3)m(N_4)} \frac{N_1N_4}{N_2^{\frac{1}{2}}N_3^{\frac{1}{2}}}N_1^{+}N_4^{+}\frac{N_1}{N_2 N_3 N_4}||\phi_1||_{X^{-1,\frac{1}{2}+}}\prod_{i=2}^{4}||\phi_i||_{X^{1,\frac{1}{2}+}}.
    \endaligned
\end{equation}
If $N \gg N_4$, then $m(N_4)=1$. Noticing $m(N_1)\leq 1$ and using \eqref{Ibasic}, we have
\begin{equation}
\aligned
    \textmd{LHS of \eqref{bound1}}&\lesssim  \frac{1}{m(N_2)m(N_3)} \frac{N_1^2 N_1^{+}N_4^{+}}{N_2^{\frac{3}{2}}N_3^{\frac{3}{2}}}||\phi_1||_{X^{-1,\frac{1}{2}+}}\prod_{i=2}^{4}||\phi_i||_{X^{1,\frac{1}{2}+}} \\
    &\lesssim    \frac{1}{N^{1-}} N_2^{0-} ||\phi_1||_{X^{-1,\frac{1}{2}+}}\prod_{i=2}^{4}||\phi_i||_{X^{1,\frac{1}{2}+}}.
    \endaligned
\end{equation}
If $N_4 \gtrsim N$, according to \eqref{Ibasic},
\begin{equation}
\aligned
    \textmd{LHS of \eqref{bound1}}&\lesssim  \frac{1}{m(N_2)m(N_3)m(N_4)} \frac{N_1^2 N_1^{+}N_4^{+}}{N_2^{\frac{3}{2}}N_3^{\frac{3}{2}}}||\phi_1||_{X^{-1,\frac{1}{2}+}}\prod_{i=2}^{4}||\phi_i||_{X^{1,\frac{1}{2}+}} \\
    &\lesssim    \frac{1}{N^{1-}} N_2^{0-} ||\phi_1||_{X^{-1,\frac{1}{2}+}}\prod_{i=2}^{4}||\phi_i||_{X^{1,\frac{1}{2}+}}.
    \endaligned
\end{equation}
Thus this subcase can be handled.\vspace{3mm}

\emph{Case III (b): $N_2 \gg N_3 \gtrsim N$. }For this subcase, we have $N_1 \sim N_2$ due to the hyperbolic plane relation. Moreover, by Bilinear estimate \eqref{bilinear2} and \eqref{Ibasic},
\begin{equation}
\aligned
    \textmd{LHS of \eqref{bound1}}&\lesssim \frac{m(N_1)}{m(N_2)m(N_3)m(N_4)} \frac{N_3N_4}{N_1^{\frac{1}{2}}N_2^{\frac{1}{2}}}N_3^{+}N_4^{+}\frac{N_1}{N_2 N_3 N_4}||\phi_1||_{X^{-1,\frac{1}{2}+}}\prod_{i=2}^{4}||\phi_i||_{X^{1,\frac{1}{2}+}} \\
    &\lesssim \frac{1}{m(N_3)m(N_4)} \frac{1}{N_2}N_3^{+}N_4^{+}||\phi_1||_{X^{-1,\frac{1}{2}+}}\prod_{i=2}^{4}||\phi_i||_{X^{1,\frac{1}{2}+}} \\
    &\lesssim    \frac{1}{N^{1-}} N_2^{0-} ||\phi_1||_{X^{-1,\frac{1}{2}+}}\prod_{i=2}^{4}||\phi_i||_{X^{1,\frac{1}{2}+}}     .
    \endaligned
\end{equation}
So this subcase can be handled as well.\vspace{3mm}

Now we estimate term $II$ in \eqref{term2}. In fact, it is easier than term $I$ in \eqref{term1}. It suffices to prove the stronger statement as follows,
\begin{equation}
     |II|=\int_0^t \int_{\Gamma_6}m_{123}(m_{456}-m_4m_5m_6)\prod_{j=1}^6 \widehat{\phi_j} \lesssim \frac{N_{max1}^{0-}}{N^{\frac{4}{3}-}}  \prod_{j=1}^6 ||\widehat{I \phi_j}||_{X^{1,\frac{1}{2}+}},
\end{equation}
where we denote $N_{max1}$ and $N_{max2}$ to be the biggest and the second biggest frequency among the $N_i$'s. Similar to the 2d case in \cite{Sgroup}, applying H\"older's inequality, \eqref{interpolate1} and \eqref{interpolate2}, we have,
\begin{equation}
    \aligned
    |II| &\lesssim    \frac{N_{max1}^{0-}}{N^{\frac{4}{3}-}} ||J^{\frac{2}{3}-}I\phi_{max1}||_{L^{\frac{30}{7}}_{t,x}}||J^{\frac{2}{3}-}I\phi_{max2}||_{L^{\frac{30}{7}}_{t,x}} \prod_{j=1}^4 ||\widehat{ \phi_j}||_{L^{\frac{15}{2}}_{t,x}}                   \\
    &\lesssim      \frac{N_{max1}^{0-}}{N^{\frac{4}{3}-}} ||I\phi_{max1}||_{X^{1,\frac{1}{2}+}}||I\phi_{max2}||_{X^{1,\frac{1}{2}+}} \prod_{j=1}^4 ||\widehat{ \phi_j}||_{X^{\frac{5}{6}+,\frac{1}{2}+}}                 \\
    &\lesssim      \frac{N_{max1}^{0-}}{N^{\frac{4}{3}-}}  \prod_{j=1}^6 ||\widehat{I \phi_j}||_{X^{1,\frac{1}{2}+}}.                 \\
    \endaligned
\end{equation}
The proof of Proposition \ref{almost} is now complete.

\begin{remark}
As mentioned in Section 3, the bilinear estimate \eqref{bilinear2} we used is stronger that bilinear estimate \eqref{bilinear3}. A stronger bilinear estimate will imply stronger local decay of the modified energy, which further implies the better low regularity result (smaller exponent $s_0$).
\end{remark}

\section{Proof of the Theorem \ref{main}}
In this section, we give the proof of our main theorem based on the I-method and the almost conservation law (Proposition \ref{almost}) established in previous sections.\vspace{3mm}

\emph{Proof of Theorem \ref{main}: }Let $u_0$ be the initial data and frequency parameter $N \gg 1$ to be decided shortly. We consider the rescaling of $u_0$ to be $u^{\lambda}_0=\frac{1}{\lambda}u_0(\frac{t}{\lambda^2},\frac{x}{\lambda},\frac{y}{\lambda})$. Choosing scaling parameter $\lambda$ such that
\begin{equation}
    \lambda \sim N^{\frac{2(1-s)}{2s-1}},
\end{equation}
we have $E(Iu^{\lambda}_0)\leq 1$. Applying Proposition \ref{almost}, for some $\delta>0$, we have
\begin{equation}
    E(Iu_{\lambda})(\delta)\lesssim  E(Iu_{\lambda})(0)+O(\frac{1}{N^{1-}}).
\end{equation}
Thus,
\begin{equation}
    E(Iu_{\lambda})(C N^{1-}\delta )\sim 1.
\end{equation}
Now we let $T$ satisfies
\begin{equation}
    T \sim \frac{C N^{1-}\delta}{\lambda^2} \sim N^{\frac{6s-5}{2s-1}-}.
\end{equation}
Hence $N$ is well-defined for all $T$ noticing that the exponent is positive when $s>\frac{5}{6}$. Now, undoing the scaling, we get that
\begin{equation}
    E(Iu_{\lambda})(T) \lesssim T^{\frac{2(1-s)}{6s-5}+}.
\end{equation}
This bound implies the desired conclusion. 

\section{Growth of higher Sobolev norm}
In this section, we focus on the Sobolev growth problem of the main equation \eqref{maineq}. There is a standard local result as follows.
\begin{proposition}\label{lwp}
(1) $($Local well-posedness$)$ Suppose $||f||_{H^1} \leq E$, then for a short time $\epsilon=\epsilon(E)\ll 1$, the equation \eqref{maineq} ($s=1$) has a unique solution $u \in X^{1,b}(-\epsilon,+\epsilon)$ with initial data $u(0)=f$, and one has
\begin{equation}
   ||u||_{X^{1,b}(-\epsilon,+\epsilon)}  \lesssim_E 1.
\end{equation}
(2) $($Propagation of regularity$)$ Moreover, if in addition $||f||_{H^2} \leq A$, then we also have
\begin{equation}
   ||u||_{X^{2,b}(-\epsilon,+\epsilon)}  \lesssim_E A.
\end{equation}
\end{proposition}
The proof of Proposition \ref{lwp} is standard and is similar to Proposition 3.1 of Deng-Germain \cite{DG}, so we omit it. The following proposition is crucial for one to establish growth of Sobolev norm result which corresponds to the `upside-down' I method.
\begin{proposition}\label{6.2}
Suppose $||u(0)||_{H^1} \leq E$ and $||\mathcal{D}u(0)||_{H^1} \leq C_1 E$ for a constant $C_1>0$,
\begin{equation}\label{growthprop}
    |E[\mathcal{D}u(T)]-E[\mathcal{D}u(0)]|\lesssim N^{-1+}+N^{o(1)}\sum_M M^{o(1)}\textmd{min}(1,N^{-1}M)||P_M\mathcal{D}u||_{L_{t,x}^{\frac{10}{3}+}}.
\end{equation}
\end{proposition}

\emph{Proof of Proposition \ref{6.2}. }Similar to the proof of Proposition \ref{almost}, direct calculations imply
\begin{equation}\label{decomp}
    \partial_t E(\mathcal{D}u)(t)=\textmd{Re}\int_{\mathcal{M}_{\lambda}}\overline{\mathcal{D}u_t}(|\mathcal{D}u|^2\mathcal{D}u-\mathcal{D}(|u|^2u)).
\end{equation}
Thus, upon integrating in time, we reduce to estimating the space-time integrals
\begin{equation}\label{Tm}
    \int_{[0,T] \times \mathcal{M}} |\mathcal{D}u|^2\overline{\mathcal{D}u} \cdot \mathcal{R}dxdt+\int_{[0,T] \times \mathcal{M}} \nabla \overline{\mathcal{D}u} \cdot \nabla \mathcal{R} dxdt,
\end{equation}
where we denote $\mathcal{R}=\mathcal{D}(|u|^2u)-|\mathcal{D}u|^2\mathcal{D}u$ for convenience.\vspace{3mm}

We start by estimating the first term in \eqref{Tm} based on a `high-low' decomposition. We denote that $u_1=P_{\leq N/10}u$ and $u_2=u-u_1=P_{> N/10}u$. Further direct calculations indicate,
\begin{equation}\label{Term1}
    \mathcal{R}=\mathcal{D}(|u|^2u-|u_1|^2u_1)-[|\mathcal{D}u|^2\mathcal{D}u-|\mathcal{D}u_1|^2\mathcal{D}u_1]-(1-\mathcal{D})(|u_1|^2u_1).
\end{equation}
For the first term in \eqref{Term1}, using H\"older and Strichartz, we can bound it by
\begin{equation}
    \aligned
    &|||\mathcal{D}u|^2\overline{\mathcal{D}u}||_{L_{t,x}^{\frac{10}{3}-}} \cdot ||\mathcal{D}(|u|^2u-|u_1|^2u_1)||_{L_{t,x}^{\frac{10}{7}+}}  \\
    &\lesssim ||\mathcal{D}u||_{L_{t,x}^{10-}}^3||\mathcal{D}u_2||_{L_{t,x}^{2+}}\big( ||\mathcal{D}u||_{L_{t,x}^{10-}}+||\mathcal{D}u_1||_{L_{t,x}^{10-}} \big)^2 \lesssim N^{-2}.
    \endaligned
\end{equation}
The bound for the second term in \eqref{Term1} holds similarly. For the last term in \eqref{Term1},
\begin{equation}
    \aligned
   ||(1-\mathcal{D}) (|u_1|^2u_1)||_{L_{t,x}^{\frac{10}{7}+}}  &\lesssim  N^{-1}||\nabla (|u_1|^2u_1)||_{L_{t,x}^{\frac{10}{7}+}}           \\
    &\lesssim  N^{-1}||u_1||_{L_{t,x}^{10-}}||\nabla  u_1||_{L_{t,x}^{2+}} \lesssim N^{-1}.
    \endaligned
\end{equation}
Gathering the above estimates, we find that
\begin{equation}
  \textmd{First term in \eqref{Tm}}   \lesssim N^{-1+}.
\end{equation}

Now we turn to the second term in \eqref{Tm}. Further direct calculations indicate,
\begin{equation}
   \nabla \mathcal{R}=[\mathcal{D}(|u|^2\nabla u)-|\mathcal{D}u|^2\nabla (\mathcal{D} u)]+[\mathcal{D}(|u|^2\nabla \bar{u})-|\mathcal{D}u|^2\nabla (\overline{\mathcal{D} u})].
\end{equation}
Since the other term is similar, we only need to consider the first term which can be decomposed as
\begin{equation}\label{Term2}
   \mathcal{D}\big( (|u|^2-|u_1|^2)\nabla u \big)-[|\mathcal{D}u|^2-|\mathcal{D} u_1|^2]\nabla(\mathcal{D}u)+\mathcal{D}(|u_1|^2\nabla u-|u_1|^2 \nabla \mathcal{D}u).
\end{equation}
For the first term in \eqref{Term2}, similar as in \cite{DG}, the following estimate holds.
\begin{equation}
    \sum_K K^{o(1)} ||\mathcal{D} P_K (|u|^2-|u_1|^2)||_{L_{t,x}^{\frac{5}{2}}} \lesssim N^{o(1)}\big( \sum_{M\geq N/10}|| \mathcal{D}P_M u||_{L_{t,x}^{\frac{10}{3}+}}+\sum_M M^{o(1)}\textmd{min}(1,N^{-1}M)||P_M \mathcal{D} u||_{L_{t,x}^{\frac{10}{3}+}} \big).
\end{equation}
Then the contribution corresponding to the first term in \eqref{Term2} can be decomposed as
\begin{equation}
    \int_{[0,T] \times \mathcal{M}} \nabla \overline{\mathcal{D}u} \cdot \mathcal{D} \big((|u|^2-|u_1|^2)\nabla u \big) dxdt=\sum_{K}\sum_{B}\int_{[0,T] \times \mathcal{M}} \nabla P_{10B} \overline{\mathcal{D}u} \cdot \mathcal{D} \big(P_K(|u|^2-|u_1|^2)\nabla P_B u \big) dxdt,
\end{equation}
where for fixed $K, B$ runs over some partition into cubes of size $K$. By orthogonality, this is bounded by
\begin{equation}
    \aligned
    & \sum_{K}\sum_{B}||\partial_i P_{10B}\overline{\mathcal{D}u}||_{L_{t,x}^{\frac{10}{3}}}||\mathcal{D} P_K (|u|^2-|u_1|^2)||_{L_{t,x}^{\frac{5}{2}}}||\mathcal{D} \partial_i P_B u||_{L_{t,x}^{\frac{10}{3}}}      \\
    &\lesssim \sum_K K^{o(1)} ||\mathcal{D} P_K (|u|^2-|u_1|^2)||_{L_{t,x}^{\frac{5}{2}}}\big(\sum_{B}||\partial_i P_{10B}\overline{\mathcal{D}u}||^2_{X^{0,b}} \big)^{\frac{1}{2}} \big(\sum_{B}||\mathcal{D} \partial_i P_B u||^2_{X^{0,b}} \big)^{\frac{1}{2}} \\
    &\lesssim N^{o(1)} \sum_M M^{o(1)}\textmd{min}(1,N^{-1}M)||P_M \mathcal{D} u||_{L_{t,x}^{\frac{10}{3}+}}.
    \endaligned
\end{equation}
The bound for the second term in \eqref{Term2} holds similarly so we omit it. For the last term in \eqref{Term2}, similar as in \cite{DG}, the following estimate holds.
\begin{equation}
   ||\mathcal{D}(P_K(|u_1|^2) \cdot \partial_i P_B u )-P_K(|u_1|^2) \cdot \partial_i \mathcal{D} P_B u||_{L^{\frac{10}{7}}_{t,x}} \lesssim \textmd{min}(1,KN^{-1})||\mathcal{D}P_K |u_1|^2||_{L^{\frac{5}{2}}_{t,x}}||\mathcal{D}\nabla P_B u||_{L^{\frac{10}{3}}_{t,x}},
\end{equation}
where we did the decomposition as before. After summing in $K$ and $B$ and using orthogonality, this gives that
\begin{equation}
    \aligned
& \big| \int_{[0,T] \times \mathcal{M}} \partial_i \overline{\mathcal{D}u} \cdot [\mathcal{D}(|u_1|^2 \partial_i u)-|u_1|^2 \partial_t \mathcal{D}u]dxdt \big|   \\
&\lesssim  \sum_{K}\sum_{B} \textmd{min}(1,N^{-1}K)  ||\mathcal{D} P_K (|u_1|^2)||_{L_{t,x}^{\frac{5}{2}}}||\nabla P_{10B}\overline{\mathcal{D}u}||_{L_{t,x}^{\frac{10}{3}}}||\mathcal{D} \nabla P_B u||_{L_{t,x}^{\frac{10}{3}}}   \\
&\lesssim  \sum_{K} \textmd{min}(1,N^{-1}K) ||\mathcal{D} P_K (|u_1|^2)||_{L_{t,x}^{\frac{5}{2}}}\big(\sum_{B}||\nabla P_{10B}\overline{\mathcal{D}u}||^2_{X^{0,b}} \big)^{\frac{1}{2}} \big(\sum_{B}||\mathcal{D} \nabla P_B u||^2_{X^{0,b}} \big)^{\frac{1}{2}}     \\
&\lesssim  \sum_{K}  \textmd{min}(1,N^{-1}K) K^{o(1)}\sum_{M\lesssim N} \textmd{min}(1,K^{-1}M)||P_M \mathcal{D} u||_{L_{t,x}^{\frac{10}{3}+}}    \\
&\lesssim     N^{o(1)} \sum_M M^{o(1)}\textmd{min}(1,N^{-1}M)||P_M \mathcal{D} u||_{L_{t,x}^{\frac{10}{3}+}}.   \\
    \endaligned
\end{equation}
This finishes the proof of Proposition \ref{growthprop}.\vspace{3mm}

At last, we give the proof of Theorem \ref{main2} as follows.\vspace{3mm}

\emph{Proof of Theorem \ref{main2}. }We choose $N=A=||u(0)||_{H^2}$. According to Sobolev embedding,
\begin{equation}
    E[\mathcal{D}u(0)] \leq C_0E,
\end{equation}
where $C_0$ is a constant. By Strichartz estimate and Proposition \ref{growthprop},
\begin{equation}
||P_M\mathcal{D}u||_{L_{t,x}^{\frac{10}{3}+}}    \lesssim M^{-1+}.
\end{equation}
As long as $T$ satisfies
\begin{equation}
\sup\limits_{0\leq t \leq T}E[\mathcal{D}u(t)] \leq 2C_0 E,
\end{equation}
using Proposition \ref{growthprop}, we have
\begin{equation}
  \sup\limits_{0\leq t \leq T}E[\mathcal{D}u(t)]-C_0 E\lesssim N^{-1+}+N^{o(1)}\sum_M M^{o(1)}\textmd{min}(1,N^{-1}M)TM^{-1+} \lesssim N^{-1+}T.
\end{equation}
The conclusion now follows from a bootstrap argument, up to $T \sim A^{1-}$.\vspace{5mm}

\noindent \textbf{Acknowledgments.} We thank Prof. Changxing Miao for some beneficial suggestions and we appreciate Chenjie Fan for some helpful discussions on their bilinear decoupling result \cite{Fan}. Zehua Zhao is partially supported by University of Maryland (postdoctoral research support).  Jiqiang Zheng was partially supported by the NSFC under grants 11771041, 11831004, 11901041. Zehua Zhao was also a guest of Institute of Applied Physics and Computational Mathematics during the writing of this paper and he appreciated for the kind host.\vspace{5mm}

 \bibliographystyle{amsplain}
 \vspace{5mm}

\noindent \author{Zehua Zhao}

\noindent \address{Department of Mathematics, University of Maryland}\\
{William E. Kirwan Hall, 4176 Campus Dr.
College Park, MD 20742-4015.}

\noindent \email{zzh@umd.edu}\\

\noindent \author{Jiqiang Zheng}

\noindent \address{Institute of Applied Physics and Computational Mathematics}\\
{Beijing, 100088, P.R.China.}

\noindent \email{zhengjiqiang@gmail.com}\\

\end{document}